# Is Global Asymptotic Stability Necessarily Uniform for Time-Invariant Time-Delay Systems?

**Iasson Karafyllis[*], Pierdomenico Pepe[**], Antoine Chaillet[***] and Yuan Wang[****]**

[*]Dept. of Mathematics, National Technical University of Athens, Zografou Campus, 15780, Athens, Greece.
emails: iasonkar@central.ntua.gr , iasonkaraf@gmail.com

[**]Dept. of Information Engineering, Computer Science, and Mathematics, University of L'Aquila, 67100, L'Aquila, Italy. email: pierdomenico.pepe@univaq.it

[***]L2S-CentraleSupélec and Université Paris-Saclay, Institut Universitaire de France, France. email: antoine.chaillet@centralesupelec.fr

[****]Dept. of Mathematical Sciences, Florida Atlantic University, Boca Raton, FL 33431, U.S.A. email: ywang@fau.edu

## Abstract

For time-invariant finite-dimensional systems, it is known that global asymptotic stability (GAS) is equivalent to uniform global asymptotic stability (UGAS), in which the decay rate and transient overshoot of solutions are requested to be uniform on bounded sets of initial states. This paper investigates this relationship for time-invariant delay systems. We show that UGAS and GAS are equivalent for this class of systems under the assumption of robust forward completeness, i.e. under the assumption that the reachable set from any bounded set of initial states on any finite time horizon is bounded. We also show that, if the state space is a space in a particular family of Sobolev or Hölder spaces, then GAS is equivalent to UGAS and that robust forward completeness holds. Based on these equivalences, we provide a novel Lyapunov characterization of GAS (and UGAS) in the aforementioned spaces.

**Keywords:** Delay Systems, Lyapunov functionals, Stability, Global Asymptotic Stability.

## 1. Introduction

For time-invariant, finite-dimensional systems described by ordinary differential equations, GAS is traditionally defined as the combination of Lyapunov stability and global convergence of solutions to the origin. An alternative way to state it is through a $KL$ bound on the solutions' norm. This alternative description is seemingly more demanding than merely stability and global attractiveness as it additionally imposes that the convergence rate and the transient overshoot of solutions are uniform over bounded sets of initial states, thus leading to the notion of Uniform Global Asymptotic Stability (UGAS). This extra conservatism turns out to be only apparent: it is well known that, for such systems, GAS and UGAS are actually equivalent properties; see [24,29].

The importance of this uniformity is twofold. First, from a practical perspective, it rules out the possibility of having an arbitrarily slow convergence of solutions to the origin or an arbitrarily large transient overshoot when initial states are confined to a bounded set. Second, it constitutes a key requirement for the construction of Lyapunov functions and is at the basis of important stability properties for systems with inputs such as Input-to-State Stability [25] and Input-to-Output Stability [27]. To that respect, it is worth mentioning that this uniformity no longer comes for free when considering output stability properties [21,14].

Another important feature of time-invariant finite-dimensional systems is that the existence of their solutions for all positive times (forward completeness) ensures a bounded reachable set over any finite time horizon from every bounded set of initial conditions [18,26]. In other words, starting from a bounded set of initial states, the solutions of a time-invariant finite-dimensional system remain bounded over a finite time horizon. In the literature, this property is referred to as either Robust Forward Completeness (RFC) [13] or bounded reachable sets property [20] and plays a crucial role in the Lyapunov characterization of forward completeness [1].

For general infinite-dimensional systems, the equivalence between GAS and UGAS and between forward completeness and RFC is far more delicate. In particular, an example is given in [19] of an infinite-dimensional system which is forward complete yet not RFC. Nevertheless, it is still an open question whether such equivalences hold when considering only time-delay systems.

Partial answers do exist though. For systems described by neutral functional differential equations, the equivalence between GAS and UGAS does not hold, even in the linear time-invariant case: see Lemma 1.1 and Example 1.6 in [10]. For time-delay systems, the relationship between GAS and UGAS was recently discussed in [23]. Interestingly, as far as local properties are concerned, asymptotic stability is indeed equivalent to uniform asymptotic stability for time-invariant delay systems [9, Lemma 1.1, p. 131]. This local result was actually proved five decades ago for globally Lipschitz time-delay systems (see Condition 4 on page 128, Definition 28.1 on page 131, Definition 30.2 on page 146, and pages 150-151 in [17]). In [15, Theorem 6.3.1, p. 73], it is stated that GAS is equivalent to UGAS for periodic delay systems provided that the function describing the dynamics is Lipschitz on bounded subsets of $C^0$. The proof is not provided in [15] and [17] is quoted for. But, as explained above, the results and the proofs provided in [17] do not show the equivalence between GAS and UGAS for the considered class of systems as given in [15].

Thus, it is not known whether GAS and UGAS are equivalent for delay systems. Similarly, it is not known whether forward completeness is equivalent to RFC for such systems, the consequences of which are discussed in [20]. Although not fully solving them, the present note shows that these two open questions are related. More specifically, we show that under the assumption of RFC, GAS and UGAS are indeed equivalent properties for time-invariant delay systems (Theorem 1). Since RFC holds automatically for globally Lipschitz delay systems and can often be established using Lyapunov techniques, our result constitutes a significant generalization of the result in [17]. The second contribution of this note is to show that the answer to both open questions depends crucially on the selection of the state space. More precisely, we show that if the considered state space is the Sobolev space $W^{1,p}\left([-r,0]\right)$ with $p \in (1,+\infty]$ then, under a forward completeness assumption, GAS is indeed equivalent to UGAS and RFC holds for time-invariant delay systems (Theorems 2 and 3). Furthermore, we show that if the considered state space is the Hölder space $C^{0,1-1/p}\left([-r,0]\right)$ with $p \in (1,+\infty]$ then, under the usual forward completeness assumption, GAS is indeed equivalent to UGAS and RFC holds, just like in the finite-dimensional case (Theorems 2 and 3). Here, it should be emphasized that Sobolev spaces have been used as state spaces in the literature for neutral delay systems: see [6,16]. The third contribution of the paper exploits this equivalence to propose a novel Lyapunov characterization of GAS (hence, UGAS) for time-invariant delay systems when treating the state space as one of the aforementioned spaces (Theorem 5 for $W^{1,p}\left([-r,0]\right)$ or $C^{0,1-1/p}\left([-r,0]\right)$ with $p \in (1,+\infty]$ and Theorem 6 for $W^{1,p}\left([-r,0]\right)$ with $p \in (1,+\infty)$ ).

The consequences of the obtained results to control theory are important. The equivalence of GAS and UGAS allows the control designer to use tools for feedback design that can prove global attractivity but not global uniform attractivity and still argue that UGAS holds. Such a tool is the extension of LaSalle's theorem in the case of delay systems (see for instance [7]), which has not been used so far for feedback control design in delay systems because it cannot guarantee uniform attractivity. Another important issue for control theory is robustness to various external inputs (disturbances). It has been shown that robustness to persistent external inputs is a consequence of uniform stability notions (see the discussion on page 162 in [17] for the case of delay systems as well as the discussion in [28] for the finite-dimensional case). The results of the present work allow the control designer to be sure that a feedback law which induces GAS for the closed-loop system will also present robustness properties with respect to various persistent external inputs.

The structure of the paper is as follows. In Section 2 we present all notions used in the paper and all main results. Section 3 provides all proofs of the main results. Finally, Section 4 concludes this paper by listing a series of open related questions.

**Notation.** Throughout this paper, we adopt the following notation.

* $\Re_+ := [0,+\infty)$. For a vector $x \in \Re^n$ we denote by $|x|$ its usual Euclidean norm.

* By $K$ we denote the set of increasing and continuous functions $\rho : \Re_+ \to \Re_+$ with $\rho(0)=0$. We say that a function $\rho \in K$ is of class $K_\infty$ if $\lim\limits_{s\to+\infty} \rho(s) = +\infty$ and $\rho(0)=0$. By $KL$ we denote the set of functions $\sigma : \Re_+ \times \Re_+ \to \Re_+$ with: (i) for each $t \geq 0$ the mapping $\sigma(\cdot,t)$ is of class $K$; (ii) for each $s \geq 0$, the mapping $\sigma(s,\cdot)$ is non-increasing with $\lim\limits_{t\to+\infty} \sigma(s,t) = 0$.

* Let $I \subseteq \Re$ be a non-empty interval and let a non-empty set $\Omega \subseteq \Re^n$. By $C^0(I\,;\Omega)$, we denote the class of continuous functions on $I$, which take values in $\Omega$. When the interval $I \subseteq \Re$ is compact, $C^0(I\,;\Omega)$ is a normed linear space with norm $\|x\|_\infty = \sup\limits_{s\in I}\left(|x(s)|\right) = \max\limits_{s\in I}\left(|x(s)|\right)$. When $I=[-r,0]$ with $r>0$ and $a \in (0,1]$ is a constant, we define the Hölder space $C^{0,a}([-r,0]) = \left\{ x \in C^0([-r,0]\,;\Re^n) : \sup\limits_{t,s\in[-r,0],t\neq s}\left(\frac{|x(t)-x(s)|}{|t-s|^a}\right) < +\infty \right\}$, i.e., the space of Hölder continuous functions of exponent $a \in (0,1]$. The Hölder space $X = C^{0,a}([-r,0])$ with $a \in (0,1]$ is a normed linear space with norm $\|x\|_X = \max\left( \|x\|_\infty , \sup\limits_{t,s\in[-r,0],t\neq s}\left(\frac{|x(t)-x(s)|}{|t-s|^a}\right)\right)$ for all $x \in C^{0,a}([-r,0])$ and by virtue of the Arzela-Ascoli theorem, for every $R>0$ the set $B = \left\{ x \in C^{0,a}([-r,0]) : \|x\|_X \leq R \right\}$ is compact in the topology of $C^0\left([-r,0];\Re^n\right)$.

* Let $I=(a,b)$ be a non-empty open interval. By $L^p\left(I;\Re^n\right)$ with $p\in[1,+\infty)$, we denote the normed linear space of equivalence classes of Lebesgue measurable functions $x:I\to\Re^n$ with $\int_a^b |x(s)|^p\,ds<+\infty$ and norm $\|x\|_p=\left(\int_a^b |x(s)|^p\,ds\right)^{1/p}$ for each $x\in L^p\left(I;\Re^n\right)$. By $L^\infty\left(I;\Re^n\right)$ we denote the normed linear space of equivalence classes of Lebesgue measurable functions $x:I\to\Re^n$ with $\sup_{a<s<b}\left(|x(s)|\right)<+\infty$ (where $\sup_{a<s<b}\left(|x(s)|\right)$ denotes the essential supremum) and norm $\|x\|_\infty=\sup_{a<s<b}\left(|x(s)|\right)$ for each $x\in L^\infty\left(I;\Re^n\right)$.

* Let $r>0$ be a given constant. We identify the Sobolev space $W^{1,p}\left([-r,0]\right)$ for $p\in[1,+\infty]$ with the normed linear space of absolutely continuous functions $x:[-r,0]\to\Re^n$ with derivative $\dot{x}$ in $L^p\left((-r,0);\Re^n\right)$. For the Sobolev space $X=W^{1,p}\left([-r,0]\right)$ we use the norm $\|x\|_X=\|x\|_\infty+\|\dot{x}\|_p$ for each $x\in W^{1,p}\left([-r,0]\right)$, which (by virtue of Theorem 8.8 on pages 212-213 in [2]) is an equivalent norm to the norm $\|x\|_p+\|\dot{x}\|_p$. Notice that $W^{1,+\infty}\left([-r,0]\right)=C^{0,1}\left([-r,0]\right)$. By virtue of Hölder's inequality, it follows that $W^{1,p}\left([-r,0]\right)\subseteq C^{0,1-1/p}\left([-r,0]\right)$ for all $p\in(1,+\infty]$. Moreover, by virtue of the Arzela-Ascoli theorem and Theorem 8.8 on pages 212-213 in [2], it follows that for every $R>0$ the bounded set $B=\left\{\,x\in W^{1,p}([-r,0]):\|x\|_\infty+\|\dot{x}\|_p\le R\,\right\}$ with $p\in(1,+\infty]$ has a closure $\bar{B}$ in $C^0\left([-r,0];\Re^n\right)$ which is compact in the topology of $C^0\left([-r,0];\Re^n\right)$ and satisfies

$$\bar{B}\subseteq\left\{\,x\in C^{0,1-1/p}([-r,0]):\max\left(\|x\|_\infty,\ \sup_{t,s\in[-r,0],t\neq s}\left(\frac{|x(t)-x(s)|}{|t-s|^{1-1/p}}\right)\right)\le R\,\right\}.$$

## 2. Main Results

### <u>2.1. Background and definitions</u>

In this work we focus on time-invariant delay systems of the form

$$\dot{x}(t)=f(x_t) \tag{1}$$

where $x(t)\in\Re^n$, $x_t\in C^0\left([-r,0];\Re^n\right)$ with $r>0$ being a constant and $(x_t)(s)=x(t+s)$ for all $s\in[-r,0]$, $f:C^0\left([-r,0];\Re^n\right)\to\Re^n$ with $f(0)=0$ is Lipschitz on bounded sets of $C^0\left([-r,0];\Re^n\right)$, i.e., there exists a non-decreasing function $L:\Re_+\to\Re_+$ such that for every $R\ge 0$ the following inequality holds

$$|f(x)-f(y)|\le L(R)\|x-y\|_\infty,\ \text{for all } x,y\in C^0\left([-r,0];\Re^n\right) \text{ with } \|x\|_\infty\le R,\|y\|_\infty\le R \tag{2}$$

Let $\phi(t,x_0)\in C^0\left([-r,0];\Re^n\right)$ denote $x_t$ where $x(t)\in\Re^n$ is the solution of (1) with initial condition $x_0\in C^0\left([-r,0];\Re^n\right)$.

The following properties have been used extensively in the literature of stability for delay systems.

**(LS)** *Lyapunov Stability: For every* $\varepsilon>0$ *there exists* $\delta(\varepsilon)>0$ *so that* $\sup\left\{\,\|\phi(t,x_0)\|_\infty : t\geq 0, x\in C^0\left([-r,0];\Re^n\right), \|x_0\|_\infty\leq\delta(\varepsilon)\,\right\}\leq\varepsilon$.

**(GA)** *Global Attractivity: For every* $x_0\in C^0\left([-r,0];\Re^n\right)$ *it holds that* $\lim\limits_{t\to+\infty}\left(\|\phi(t,x_0)\|_\infty\right)=0$.

**(UGA)** *Uniform Global Attractivity: For every* $\varepsilon,\rho>0$ *there exists* $T(\varepsilon,\rho)>0$ *so that* $\sup\left\{\,\|\phi(t,x_0)\|_\infty : t\geq T(\varepsilon,\rho), x_0\in C^0\left([-r,0];\Re^n\right), \|x_0\|_\infty\leq\rho\,\right\}\leq\varepsilon$.

**(LagS)** *Lagrange Stability: For every* $\rho>0$ *it holds that* $\sup\left\{\,\|\phi(t,x_0)\|_\infty : t\geq 0, x_0\in C^0\left([-r,0];\Re^n\right), \|x_0\|_\infty\leq\rho\,\right\}<+\infty$.

**(RFC)** *Robust Forward Completeness: For every* $\rho,T>0$ *it holds that* $\sup\left\{\,\|\phi(t,x_0)\|_\infty : t\in[0,T], x_0\in C^0\left([-r,0];\Re^n\right), \|x_0\|_\infty\leq\rho\,\right\}<+\infty$.

**(GAS)** *Global Asymptotic Stability: Both properties (LS) and (GA) hold.*

**(UGAS)** *Uniform Global Asymptotic Stability: There exists* $\sigma\in KL$ *such that the estimate* $\|\phi(t,x_0)\|_\infty\leq\sigma\left(\,\|x_0\|_\infty, t\,\right)$ *holds for all* $t\geq 0$ *and all* $x_0\in C^0\left([-r,0];\Re^n\right)$.

Lyapunov Stability (LS) is a purely local property and imposes that solutions remain arbitrarily close to the origin provided that the norm of the initial segment is sufficiently small. Property (GA) requires that all solutions eventually converge to the origin. Property (UGA) additionally requires that the convergence rate at which solutions converge is uniform for bounded sets of initial states. Lagrange stability (LagS) can be interpreted as solutions' boundedness. Robust Forward Completeness requires not only existence of solutions for all forward times, but also that their magnitude is bounded over any compact time interval and for initial states in any bounded set. Finally, (UGAS) employs the classical $KL$ formalism and readily implies both (LS) and (UGA). Some of the above properties are related. For example, it is well-known (see Theorem 2.2 on page 62 in [13]) that the implication (LS) $\wedge$ (UGA) $\wedge$ (LagS) $\Leftrightarrow$ (UGAS) holds. Moreover, Lemma 2.1 on page 58 in [13] shows that the implication (UGA) $\wedge$ (RFC) $\Leftrightarrow$ (UGAS) holds. Some of the above properties are stronger than others; for example, the implications (UGA) $\Rightarrow$ (GA), (LagS) $\Rightarrow$ (RFC) hold trivially.

However, to the best of our knowledge, it is not known whether the implication (GAS) $\Rightarrow$ (UGAS) is true or not for delay systems. This implication is true for delay-free (finite-dimensional) systems (see [24,29]). Another important property that is valid for delay-free systems (see [18,26]) is the property that simple forward completeness (i.e., global existence of solutions for arbitrary initial condition) implies Robust Forward Completeness (RFC). Whether or not such equivalence also holds for time-delay systems is another open question.

*2.2. Uniformity under RFC*

The first contribution of this paper is to show that, provided that RFC holds, (GAS) is equivalent to (UGAS).

**Theorem 1:** *The following implications hold for system (1):*

$$\text{(GAS)} \wedge \text{(RFC)} \Leftrightarrow \text{(UGAS)}$$

In other words, provided that RFC holds, the combination of Lyapunov stability and global attractivity does ensure Uniform Global Asymptotic Stability, just like in the finite-dimensional case. It is worth stressing that:

(i) RFC holds automatically for important classes of delay systems such as globally Lipschitz delay systems. This reminds the results proved in [17] for globally Lipschitz systems (see Condition 4, p. 128, Definition 28.1, p. 131, Definition 30.2, p. 146, and pages 150, 151 in [17]).

(ii) RFC can often be established by using Lyapunov-like functionals (see [1] for the finite-dimensional case and [13] for the time-delay case). For example, the existence of a functional $U : C^0\left([-r,0];\Re^n\right) \to \Re_+$ which is Lipschitz on bounded sets of $C^0\left([-r,0];\Re^n\right)$ and for which there exist a function $a \in K_\infty$ and a constant $\mu \geq 0$ such that the inequalities $U(x) \geq a\left(|x(0)|\right)$ and $\limsup_{h \to 0^+}\left(h^{-1}\left(U\left(P_h(x)\right) - U(x)\right)\right) \leq \mu U(x)$ hold for all $x \in C^0\left([-r,0];\Re^n\right)$, where $\left(P_h(x)\right)(s) = x(s+h)$ for $h \geq 0$ and $s \in [-r,-h]$ and $\left(P_h(x)\right)(s) = x(0) + (s+h) f(x)$ for $h > 0$ and $s \in (-h,0]$, is sufficient to guarantee RFC.

Since (GA) obviously implies forward completeness, the RFC requirement in Theorem 1 could be removed if one could establish that forward completeness implies RFC for time-delay systems (as it holds in finite dimension). To date, this crucial question remains open.

The fact that Lagrange stability trivially ensures robust forward completeness, gives us the following corollary. Notice that, differently from [15, Theorem 6.3.1, p. 73], here the (LagS) property is invoked. Notice also that (LS) $\wedge$ (LagS) is equivalent to (LS) with the additional requirement that $\delta(\varepsilon)$ for which $\sup\left\{ \|\phi(t,x_0)\|_\infty : t \geq 0, x \in C^0\left([-r,0];\Re^n\right), \|x_0\|_\infty \leq \delta(\varepsilon) \right\} \leq \varepsilon$ can be chosen arbitrarily large for sufficiently large $\varepsilon > 0$.

**Corollary 1:** *The following implications hold for system (1):*

$$\text{(GAS)} \wedge \text{(LagS)} \Leftrightarrow \text{(UGAS)}$$

*2.3. Uniformity in Sobolev and Hölder spaces*

The second contribution of this paper is to show that, when working in a Sobolev space $W^{1,p}\left([-r,0]\right)$ with $p \in (1,+\infty]$ or a Hölder space $C^{0,q}\left([-r,0]\right)$ with $q \in (0,1]$, the RFC requirement can be removed from the above implications. It is a fact that if the initial condition $x_0 \in C^0\left([-r,0];\Re^n\right)$ is of class $W^{1,p}\left([-r,0]\right)$ for some $p \in [1,+\infty]$ then the solution $\phi(t,x_0)$ is of class $W^{1,p}\left([-r,0]\right)$ whenever it exists. Similarly, if the initial condition $x_0 \in C^0\left([-r,0];\Re^n\right)$ is of class $C^{0,q}([-r,0])$ for some $q\in(0,1]$ then the solution $\phi(t,x_0)$ is of class $C^{0,q}([-r,0])$. These facts have been utilized in some works on delay systems (see for example [11,12,22] with $p=+\infty$ or $q=1$). Therefore, instead of considering the state space to be $C^0([-r,0];\Re^n)$ we may consider as state space the normed linear space $X=W^{1,p}([-r,0])$ for some $p\in(1,+\infty]$ with $\|x\|_X=\|x\|_\infty+\|\dot{x}\|_p$ for each $x\in W^{1,p}([-r,0])$ or the normed linear space $X=C^{0,q}([-r,0])$ for some $q\in(0,1]$ with $\|x\|_X=\max\left(\|x\|_\infty,\ \sup_{t,s\in[-r,0],t\neq s}\left(\frac{|x(t)-x(s)|}{|t-s|^q}\right)\right)$ for each $x\in C^{0,q}([-r,0])$.

The change of the state space requires updating of the properties listed above by replacing the sup norm by the Sobolev norm or the Hölder norm. This leads to the following counterparts in which $X$ may denote either the Sobolev space $W^{1,p}([-r,0])$ for some $p\in(1,+\infty]$ or the Hölder space $C^{0,q}([-r,0])$ for some $q\in(0,1]$.

**(LS-X)** *Lyapunov Stability: For every* $\varepsilon>0$ *there exists* $\delta(\varepsilon)>0$ *so that* $\sup\left\{\|\phi(t,x_0)\|_X : t\geq 0, x_0\in X, \|x_0\|_X\leq\delta(\varepsilon)\right\}\leq\varepsilon$.

**(GA-X)** *Global Attractivity: For every* $x_0\in X$ *it holds that* $\lim_{t\to+\infty}\left(\|\phi(t,x_0)\|_X\right)=0$.

**(UGA-X)** *Uniform Global Attractivity: For every* $\varepsilon,\rho>0$ *there exists* $T(\varepsilon,\rho)>0$ *so that* $\sup\left\{\|\phi(t,x_0)\|_X : t\geq T(\varepsilon,\rho), x_0\in X, \|x_0\|_X\leq\rho\right\}\leq\varepsilon$.

**(LagS-X)** *Lagrange Stability: For every* $\rho>0$ *it holds that* $\sup\left\{\|\phi(t,x_0)\|_X : t\geq 0, x_0\in X, \|x_0\|_X\leq\rho\right\}<+\infty$.

**(RFC-X)** *Robust Forward Completeness: For every* $\rho,T>0$ *it holds that* $\sup\left\{\|\phi(t,x_0)\|_X : t\in[0,T], x_0\in X, \|x_0\|_X\leq\rho\right\}<+\infty$.

**(GAS-X)** *Global Asymptotic Stability: Both properties (LS-X) and (GA-X) hold.*

**(UGAS-X)** *Uniform Global Asymptotic Stability: There exists* $\sigma\in KL$ *such that the estimate* $\|\phi(t,x_0)\|_X\leq\sigma\left(\|x_0\|_X, t\right)$ *holds for all* $t\geq 0$ *and all* $x_0\in X$.

With these definitions at hand, we are able to show that simple forward completeness (i.e., global existence of solutions for arbitrary initial condition) implies RFC in particular Sobolev or Hölder spaces.

**Theorem 2:** *Let* $p\in(1,+\infty]$ *be given. Suppose that (1) is forward complete, in the sense that for every* $x_0\in C^{0,1-1/p}([-r,0])$ *the solution* $x(t)\in\Re^n$ *of (1) with initial condition* $x_0$ *exists for all* $t\geq 0$. *Then (1) with state space* $X=W^{1,p}([-r,0])$ *or state space* $X=C^{0,1-1/p}([-r,0])$ *is Robustly Forward Complete, i.e., Property (RFC-X) holds.*

Moreover, in this case we are able to show that Lyapunov stability combined with (non-uniform) global attractivity implies Uniform Global Asymptotic Stability.

**Theorem 3:** *Let* $X$ *denote the Sobolev space* $W^{1,p}\left([-r,0]\right)$ *or the Hölder space* $C^{0,1-1/p}\left([-r,0]\right)$ *for some* $p\in(1,+\infty]$. *Suppose that (1) is forward complete, in the sense that for every* $x_0\in C^{0,1-1/p}\left([-r,0]\right)$ *the solution* $x(t)\in\Re^n$ *of (1) with initial condition* $x_0$ *exists for all* $t\geq 0$. *Then the following implications hold for system (1):*

$$\text{(GAS-X)} \Leftrightarrow \text{(UGAS-X)}$$

Clearly, when $X=C^{0,1-1/p}\left([-r,0]\right)$, the assumption in Theorem 3 that for every $x_0\in C^{0,1-1/p}\left([-r,0]\right)$ the solution $x(t)\in\Re^n$ of (1) with initial condition $x_0\in C^{0,1-1/p}\left([-r,0]\right)$ exists for all $t\geq 0$ is a redundant assumption (since both (GAS-X) and (UGAS-X) imply this property). Thus, when working with the Hölder spaces $X=C^{0,1-1/p}\left([-r,0]\right)$ with $p\in(1,+\infty]$, the combination of Lyapunov stability and global attractivity is equivalent to (UGAS), just like for finite-dimensional systems.

The stability properties of system (1) viewed in different state spaces are related. The following theorem uses the following stability notion, which provides a $KL$ bound on the sup norm of the state in terms of the Sobolev norm (when $X=W^{1,p}\left([-r,0]\right)$) or the Hölder norm (when $X=C^{0,1-1/p}\left([-r,0]\right)$) of the initial condition.

**(Q-X)** *There exists* $\sigma\in KL$ *such that the estimate* $\left\|\phi(t,x_0)\right\|_\infty\leq\sigma\left(\left\|x_0\right\|_X,t\right)$ *holds for all* $t\geq 0$, $x_0\in X$.

**Theorem 4:** *Let* $X$ *denote the Sobolev space* $W^{1,p}\left([-r,0]\right)$ *or the Hölder space* $X=C^{0,1-1/p}\left([-r,0]\right)$ *for some* $p\in(1,+\infty]$. *The following implications hold for system (1):*

$$\text{(UGAS)} \Rightarrow \text{(Q-X)} \Leftrightarrow \text{(UGAS-X)}$$

A direct consequence of the inequality $\left\|\dot{x}_0\right\|_p\leq r^{\frac{1}{p}-\frac{1}{q}}\left\|\dot{x}_0\right\|_q$, which holds for any $q\geq p\geq 1$ and all $\dot{x}_0\in L^q\left((-r,0);\Re^n\right)$ (obtained by using Hölder's inequality), the fact that $\left\|x_0\right\|_X=\left\|x_0\right\|_\infty+\left\|\dot{x}_0\right\|_p$ when $X=W^{1,p}\left([-r,0]\right)$ and Theorem 4 is the following corollary.

**Corollary 2:** *Suppose that the (UGAS-X) property holds with* $X$ *being the space* $W^{1,p}\left([-r,0]\right)$ *for some* $p\in(1,+\infty]$. *Then for every* $q\in[p,+\infty]$ *the (UGAS-X) property holds with* $X$ *being the space* $W^{1,q}\left([-r,0]\right)$.

A direct consequence of the inequality $\sup\limits_{l,s\in[-r,0],l\neq s}\left(\frac{\left|x_0(l)-x_0(s)\right|}{\left|l-s\right|^{1-1/p}}\right)\leq r^{\frac{1}{p}-\frac{1}{q}}\sup\limits_{l,s\in[-r,0],l\neq s}\left(\frac{\left|x_0(l)-x_0(s)\right|}{\left|l-s\right|^{1-1/q}}\right)$, which holds for each $q\geq p>1$ and all $x_0\in C^{0,1-1/q}\left([-r,0]\right)$, the fact that

$\|x_0\|_X = \max\left( \|x_0\|_\infty , \sup_{t,s\in[-r,0],t\neq s}\left( \frac{|x_0(t)-x_0(s)|}{|t-s|^{1-1/p}} \right) \right)$ when $X = C^{0,1-1/p}([-r,0])$ and Theorem 4 is the following corollary.

**Corollary 3:** *Suppose that the (UGAS-X) property holds with $X$ being the space $C^{0,1-1/p}([-r,0])$ for some $p \in (1,+\infty]$. Then for every $q \in [p,+\infty]$ the (UGAS-X) property holds with $X$ being the space $C^{0,1-1/q}([-r,0])$.*

A direct consequence of the inequality $\sup_{l,s\in[-r,0],l\neq s}\left( \frac{|x_0(l)-x_0(s)|}{|l-s|^{1-1/p}} \right) \le \|\dot{x}_0\|_p$, which holds for all $p>1$ and all $x_0 \in W^{1,p}([-r,0])$, the facts that $\|x_0\|_X = \max\left( \|x_0\|_\infty , \sup_{t,s\in[-r,0],t\neq s}\left( \frac{|x_0(t)-x_0(s)|}{|t-s|^{1-1/p}} \right) \right)$ when $X = C^{0,1-1/p}([-r,0])$ and $\|x_0\|_X = \|x_0\|_\infty + \|\dot{x}_0\|_p$ when $X = W^{1,p}([-r,0])$ and Theorem 4 is the following corollary.

**Corollary 4:** *Suppose that the (UGAS-X) property holds with $X$ being the Hölder space $C^{0,1-1/p}([-r,0])$ for some $p \in (1,+\infty]$. Then the (UGAS-X) property holds with $X$ being the Sobolev space $W^{1,p}([-r,0])$.*

## <u>2.4. Lyapunov-Krasovskii characterization of GAS in Sobolev and Hölder spaces</u>

This section contains our third contribution. Using fundamental properties of delay systems and the converse Lyapunov theory in [13] we obtain the following Lyapunov characterization of the (UGAS-X) property when $X$ denotes either a Sobolev or a Hölder space.

**Theorem 5:** *Let $X$ denote the space $W^{1,p}([-r,0])$ or the space $C^{0,1-1/p}([-r,0])$ for some $p \in (1,+\infty]$. Property (UGAS-X) holds if and only if there exist a functional $V : X \to \Re_+$ which is Lipschitz on bounded sets of $X$ and functions $a_1, a_2 \in K_\infty$ such that the following inequalities hold for all $x \in X$ :*

$$a_1\left(\|x\|_X\right) \le V(x) \le a_2\left(\|x\|_X\right) \tag{3}$$

$$V(\phi(t,x)) \le \exp(-t)V(x) \text{, for all } t \ge 0 . \tag{4}$$

The strength of Theorem 5 lies in the fact that the constructed Lyapunov-Krasovskii functional is coercive (in the sense that it is sandwiched between two $K_\infty$ functions of the Sobolev norms or Hölder norms of the state) and that it decays exponentially fast along solutions. The problem with Theorem 5 is that we were not able to obtain a differential inequality that is equivalent to inequality (4). Indeed, (4) implies the differential inequality

$$\limsup_{t\to 0^+}\left( t^{-1}\left( V(\phi(t,x)) - V(x) \right) \right) \le -V(x) \tag{5}$$

but the differential inequality (5) does not necessarily imply (4) since we are not aware whether the mapping $t \mapsto V(\phi(t,x))$ is lower semi-continuous or not (if the mapping $t \mapsto V(\phi(t,x))$ were lower semi-continuous then an application of Lemma 2.12 on pages 77-78 in [13] would imply (4)). More specifically, since we are not aware if the mapping $t \mapsto \phi(t,x)$ is continuous in the topology of $W^{1,\infty}([-r,0])$ or the topology of $C^{0,1-1/p}([-r,0])$ for $p \in (1,+\infty]$, the fact that the Lyapunov functional $V$ is Lipschitz on bounded sets of $W^{1,\infty}([-r,0])$ or on bounded sets of $C^{0,1-1/p}([-r,0])$ with $p \in (1,+\infty]$ cannot guarantee continuity (or lower semi-continuity) for the mapping $t \mapsto V(\phi(t,x))$. However, when $p < +\infty$ the mapping $t \mapsto \phi(t,x)$ is in fact continuous in the topology of $W^{1,p}([-r,0])$. This is guaranteed by the following lemma.

**Lemma 1:** *Let $X$ denote the space $W^{1,p}([-r,0])$ for some $p \in (1,+\infty)$. Then given any $x_0 \in X$, the mapping*

$$t \mapsto \phi(t,x_0) \tag{6}$$

*is continuous on $[0,t_{\max})$ in the topology of $X$, where $t_{\max} \in (0,+\infty]$ is the maximal existence time of the solution of (1) with initial condition $x_0 \in X$.*

Therefore, when $X = W^{1,p}([-r,0])$ for some $p \in (1,+\infty)$, we are able to characterize the (GAS-X) (hence, the (UGAS-X)) property through a coercive Lyapunov-Krasovskii functional whose upper right Dini derivative leads to an exponential decay estimate.

**Theorem 6:** *Let $X$ denote the space $W^{1,p}([-r,0])$ for some $p \in (1,+\infty)$. Suppose that (1) is forward complete, in the sense that for every $x_0 \in C^{0,1-1/p}([-r,0])$ the solution $x(t) \in \Re^n$ of (1) with initial condition $x_0 \in C^{0,1-1/p}([-r,0])$ exists for all $t \geq 0$. Then the following statements are equivalent for system (1).*

***(i)*** *Property (GAS-X) holds.*

***(ii)*** *Property (UGAS-X) holds.*

***(iii)*** *There exist a functional $V : X \to \Re_+$ which is Lipschitz on bounded sets of $X$ and functions $a_1, a_2 \in K_\infty$ such that the following inequalities hold for all $x \in X$:*

$$a_1\left(\|x\|_X\right) \leq V(x) \leq a_2\left(\|x\|_X\right) \tag{7}$$

$$\limsup_{t \to 0^+}\left(t^{-1}\left(V(\phi(t,x)) - V(x)\right)\right) \leq -V(x) \tag{8}$$

***(iv)*** *There exist a functional $V : X \to \Re_+$ which is Lipschitz on bounded sets of $X$, a continuous positive definite function $Q : \Re^n \to \Re_+$ and functions $a_1, a_2 \in K_\infty$ such that the following inequalities hold for all $x \in X$:*

$$a_1\left(|x(0)|\right) \leq V(x) \leq a_2\left(\|x\|_X\right) \tag{9}$$

$$\limsup_{t \to 0^+}\left(t^{-1}\left(V(\phi(t,x)) - V(x)\right)\right) \leq -Q(x(0)) \tag{10}$$

It can be easily checked that (7)-(8) imply (9)-(10). Indeed, (7)-(8) are conditions for a coercive Lyapunov-Krasovskii functional with dissipation rate involving the whole functional itself, whereas the functional in (9)-(10) is not required to be coercive (its lower bound in (9) involves only the current value of the solution rather than the Sobolev norm of the history of the solution) and its dissipation rate merely involves the current value of the solution. As discussed in [3] and [4], the use of non-coercive Lyapunov-Krasovskii functional with point-wise dissipation rate is often more convenient to ensure GAS in practice. On the other hand, coercive Lyapunov-Krasovskii functionals with exponential decay rate constitute a valuable tool for robustness analysis.

The assumption that (1) is forward complete, in the sense that for every $x_0 \in C^{0,1-1/p}\left([-r,0]\right)$ the solution $x(t)\in\Re^n$ of (1) with initial condition $x_0 \in C^{0,1-1/p}\left([-r,0]\right)$ exists for all $t\geq 0$, is not redundant in Theorem 6. This assumption is needed only for the implication (i) $\Rightarrow$ (ii) (guaranteed by means of Theorem 3). Indeed, the property (GAS-X) does not guarantee the required forward completeness assumption when $X=W^{1,p}\left([-r,0]\right)$. Notice that each space $W^{1,p}\left([-r,0]\right)$ with $p\in(1,+\infty)$ is a strict subset of the space $C^{0,1-1/p}\left([-r,0]\right)$ and property (GAS-X) guarantees only the existence of global solutions when the initial condition is in $X=W^{1,p}\left([-r,0]\right)$.

Finally, it should be noticed that Theorems 2, 3, 4, 5 and 6 exclude the case $X=W^{1,p}\left([-r,0]\right)$ with $p=1$ (Theorem 6 also excludes the case $p=+\infty$ for the reasons that were explained above). There is a fundamental reason for this exclusion: the closure of the unit ball of $W^{1,p}\left([-r,0]\right)$ is compact in the topology of $C^0\left([-r,0];\Re^n\right)$ when $p\in(1,+\infty]$ but is not compact in the topology of $C^0\left([-r,0];\Re^n\right)$ when $p=1$ (see Theorem 8.8 on pages 212-213 in [2]). Therefore, we are not aware whether these theorems hold for the case $X=W^{1,p}\left([-r,0]\right)$ with $p=1$.

## 3. Proofs

In the following proofs, we use the following facts that hold for any $a<b$ and any $c\in(a,b)$:

- Given $p\in(1,+\infty)$, if $x:(a,b)\to\Re^n$ is a Lebesgue measurable function which is of class $L^p\left((a,c)\right)\cap L^\infty\left((c,b)\right)$ then
$$\|x\|_p \leq \left(1+(b-a)^{1/p}\right)\max\left(\left(\int_a^c |x(s)|^p\,ds\right)^{1/p},\ \sup_{s\in(c,b)}\left(|x(s)|\right)\right).$$

- Given $p\in(1,+\infty]$, if $x:(a,b)\to\Re^n$ is an essentially bounded Lebesgue measurable function then
$$\|x\|_p \leq (b-a)^{1/p}\,\|x\|_\infty .$$

- Given $p\in(1,+\infty]$, if $x:[a,b]\to\Re^n$ is a continuous function which is absolutely continuous on $[c,b]$ with $\sup_{l,s\in[a,c],l\neq s}\left(\frac{|x(l)-x(s)|}{|l-s|^{1-1/p}}\right)<+\infty$ and $\sup_{s\in(c,b)}\left(|\dot{x}(s)|\right)<+\infty$ then

$$\sup_{l,s\in[a,b],l\neq s}\left(\frac{|x(l)-x(s)|}{|l-s|^{1-1/p}}\right)\leq \sup_{l,s\in[a,c],l\neq s}\left(\frac{|x(l)-x(s)|}{|l-s|^{1-1/p}}\right)+(b-a)^{1/p}\sup_{s\in(c,b)}\left(|\dot{x}(s)|\right).$$

- Given $p\in(1,+\infty]$, if $x:[a,b]\to\Re^n$ is a continuous function which is absolutely continuous on $[a,b]$ with $\sup_{s\in(a,b)}\left(|\dot{x}(s)|\right)<+\infty$ then

$$\sup_{l,s\in[a,b],l\neq s}\left(\frac{|x(l)-x(s)|}{|l-s|^{1-1/p}}\right)\leq (b-a)^{1/p}\sup_{s\in(a,b)}\left(|\dot{x}(s)|\right).$$

**Proof of Theorem 1:** Using Lemma 2.1 on page 58 in [13] (i.e., the equivalence (UGA) $\wedge$ (RFC) $\Leftrightarrow$ (UGAS)), and Theorem 2.2 on page 62 in [13]) (i.e., the equivalence (LS) $\wedge$ (UGA) $\wedge$ (LagS) $\Leftrightarrow$ (UGAS)) it suffices to prove the implication (LS) $\wedge$ (GA) $\wedge$ (RFC) $\Rightarrow$ (UGA).

Let arbitrary $\varepsilon,\rho>0$ be given.

By virtue of Lyapunov Stability (Property (LS)), there exists $\delta>0$ so that

$$\sup\left\{ \|\phi(t,x_0)\|_\infty : t\geq 0, x_0\in C^0\left([-r,0];\Re^n\right), \|x_0\|_\infty\leq\delta \right\}\leq\varepsilon \tag{11}$$

By virtue of RFC there exists $R>0$ so that

$$\sup\left\{ \|\phi(t,x_0)\|_\infty : t\in[0,r], x_0\in C^0\left([-r,0];\Re^n\right), \|x_0\|_\infty\leq\rho \right\}\leq R \tag{12}$$

Let $M:=\max\left\{ |f(x)| : \|x\|_\infty\leq R \right\}$ and define the set

$$B:=\left\{ x\in W^{1,\infty}\left([-r,0];\Re^n\right) : \|\dot{x}\|_\infty\leq M, \|x\|_\infty\leq R \right\} \tag{13}$$

Then by virtue of the Arzela-Ascoli theorem and the fact that a function $x:[-r,0]\to\Re^n$ is in $B$ if and only if $x$ is Lipschitz with Lipschitz constant $M$ and $\|x\|_\infty\leq R$, the set $B$ is compact in the topology of $C^0\left([-r,0];\Re^n\right)$.

Definition (13), estimate (12) and (1) imply that for every $x_0\in C^0\left([-r,0];\Re^n\right)$ with $\|x_0\|_\infty\leq\rho$ it holds that $\phi(r,x_0)\in B$.

Let arbitrary $x_0\in C^0\left([-r,0];\Re^n\right)$ with $\|x_0\|_\infty\leq R$ be given. By virtue of Global Attractivity (Property (GA)), there exists $\tilde{T}(x_0)>0$ such that $\|\phi(t,x_0)\|_\infty\leq\frac{\delta}{2}$ for all $t\geq\tilde{T}(x_0)$.

By virtue of continuity of solutions with respect to the initial conditions for each $x_0\in C^0\left([-r,0];\Re^n\right)$ with $\|x_0\|_\infty\leq R$ there exists $\eta(x_0)>0$ such that $\|\phi(t,y)-\phi(t,x_0)\|_\infty\leq\frac{\delta}{2}$ for all $t \in \left[0, \tilde{T}(x_0)\right]$ and for all $y \in C^0\left([-r,0];\Re^n\right)$ with $\|y - x_0\|_\infty < \eta(x_0)$. Thus, by virtue of the triangle inequality, we get for all $y \in C^0\left([-r,0];\Re^n\right)$ with $\|y - x_0\|_\infty < \eta(x_0)$:

$$\left\|\phi\left(\tilde{T}(x_0), y\right)\right\|_\infty \leq \left\|\phi\left(\tilde{T}(x_0), y\right) - \phi\left(\tilde{T}(x_0), x_0\right)\right\|_\infty + \left\|\phi\left(\tilde{T}(x_0), x_0\right)\right\|_\infty \leq \frac{\delta}{2} + \frac{\delta}{2} \leq \delta \tag{14}$$

Consequently, by virtue of (11), (14) and the semigroup property we get for each $x_0 \in C^0\left([-r,0];\Re^n\right)$ with $\|x_0\|_\infty \leq R$:

$$\begin{aligned}
&\sup\left\{ \|\phi(t,y)\|_\infty : t \geq \tilde{T}(x_0),\, y \in C^0\left([-r,0];\Re^n\right), \|y - x_0\|_\infty < \eta(x_0) \right\} \\
&= \sup\left\{ \left\|\phi\left(t - \tilde{T}(x_0), \phi\left(\tilde{T}(x_0), y\right)\right)\right\|_\infty : t \geq \tilde{T}(x_0),\, y \in C^0\left([-r,0];\Re^n\right), \|y - x_0\|_\infty < \eta(x_0) \right\} \\
&\leq \sup\left\{ \|\phi(s,z)\|_\infty : s \geq 0,\, z \in C^0\left([-r,0];\Re^n\right), \|z\|_\infty \leq \delta \right\} \leq \varepsilon
\end{aligned} \tag{15}$$

It follows that for each $x \in B$ the set $N(x) := \left\{ y \in C^0\left([-r,0];\Re^n\right) : \|y - x\|_\infty < \eta(x) \right\}$ is an open neighborhood of $x$ in $C^0\left([-r,0];\Re^n\right)$. Consequently, the sets $N(x) \cap B$ for $x \in B$ constitute an open cover of $B$. Therefore, by compactness of the set $B$ there exists a (finite) positive integer $m$ and points $x_i \in B$, $i = 1,...,m$ such that $B \subseteq \bigcup_{i=1}^{m} N(x_i)$.

Define

$$\begin{aligned}
&T = T(\varepsilon, \rho) := r + \tau \\
&\tau := \max_{i=1,...,m}\left(\tilde{T}(x_i)\right)
\end{aligned} \tag{16}$$

Let arbitrary $y \in B$ be given. Since $B \subseteq \bigcup_{i=1}^{m} N(x_i)$, there exists $i \in \{1,...,m\}$ such that $y \in N(x_i) = \left\{ z \in C^0\left([-r,0];\Re^n\right) : \|z - x_i\|_\infty \leq \eta(x_i) \right\}$. Consequently, by virtue of definition (16) and inequality (15), we get:

$$\sup\left\{ \|\phi(t,y)\|_\infty : t \geq \tau \right\} \leq \sup\left\{ \|\phi(t,z)\|_\infty : t \geq \tilde{T}(x_i),\, z \in C^0\left([-r,0];\Re^n\right), \|z - x_i\|_\infty < \eta(x_i) \right\} \leq \varepsilon \tag{17}$$

It follows from (17) and the fact that $y \in B$ is arbitrary that the following inequality holds:

$$\sup\left\{ \|\phi(t,y)\|_\infty : t \geq \tau,\, y \in B \right\} \leq \varepsilon \tag{18}$$

The fact that for every $x_0 \in C^0\left([-r,0];\Re^n\right)$ with $\|x_0\|_\infty \leq \rho$ it holds that $\phi(r, x_0) \in B$ in conjunction with (18), definition (16) and the semigroup property give:

$$
\begin{aligned}
&\sup\left\{ \|\phi(t,x_0)\|_\infty : t \geq T\,, x_0 \in C^0\left([-r,0];\Re^n\right), \|x_0\|_\infty \leq \rho \right\} \\
&= \sup\left\{ \|\phi(t,x_0)\|_\infty : t \geq r+\tau\,, x_0 \in C^0\left([-r,0];\Re^n\right), \|x_0\|_\infty \leq \rho \right\} \\
&= \sup\left\{ \left\|\phi\left(s,\phi(r,x)\right)\right\|_\infty : s \geq \tau\,, x_0 \in C^0\left([-r,0];\Re^n\right), \|x_0\|_\infty \leq \rho \right\} \\
&\leq \sup\left\{ \|\phi(t,y)\|_\infty : t \geq \tau\,, y \in B \right\} \leq \varepsilon
\end{aligned} \tag{19}
$$

Inequality (19) directly implies Uniform Global Attractivity (Property (UGA)).

The proof is complete. $\triangleleft$

**Proof of Theorem 2:** We want to prove that for every $R,T>0$ it holds that $\sup\left\{ \|x_t\|_X : t\in[0,T], x_0 \in X\,, \|x_0\|_X \leq R \right\} < +\infty$. In order to prove this, it suffices to prove that for every $R>0$ it holds that

$$\sup\left\{ \|x_t\|_X : t\in[0,r], x_0 \in X\,, \|x_0\|_X \leq R \right\} < +\infty \tag{20}$$

To see this, notice that assuming that (20) holds and applying (20) inductively in conjunction with the semigroup property, we get

$$\sup\left\{ \|x_t\|_X : t\in[0,nr], x_0 \in X\,, \|x_0\|_X \leq R \right\} < +\infty, \text{ for } n=1,2,... \tag{21}$$

It is clear that (21) implies $\sup\left\{ \|x_t\|_X : t\in[0,T], x_0 \in X\,, \|x_0\|_X \leq R \right\} < +\infty$ for every $R,T>0$, i.e., Property (RFC-X).

Furthermore, in order to show (20), it suffices to prove that for every $R>0$ it holds that

$$\rho = \sup\left\{ \|x_t\|_\infty : t\in[0,r], x_0 \in X\,, \|x_0\|_X \leq R \right\} < +\infty \tag{22}$$

Indeed, assuming that (22) holds and using (1), (2) and the fact that $f(0)=0$, we get for $t\in[0,r]$ that

$$\sup\left\{ \|\dot{x}_t\|_p : t\in[0,r], x_0 \in X\,, \|x_0\|_X \leq R \right\} \leq \left(1+r^{1/p}\right)\max\left(R, \rho L(\rho)\right)$$
$$\text{when } X = W^{1,p}\left([-r,0]\right) \text{ with } p\in(1,+\infty)$$

$$\sup\left\{ \|\dot{x}_t\|_p : t\in[0,r], x_0 \in X\,, \|x_0\|_X \leq R \right\} \leq \max\left(R, \rho L(\rho)\right)$$
$$\text{when } X = W^{1,p}\left([-r,0]\right) \text{ with } p=+\infty$$

$$\sup\left\{ \sup_{l,s\in[-r,0], l\neq s}\left( \frac{|x_t(l)-x_t(s)|}{|l-s|^{1-1/p}} \right) : t\in[0,r], x_0 \in X\,, \|x_0\|_X \leq R \right\} \leq R + r^{1/p}\rho L(\rho)$$
$$\text{when } X = C^{0,1-1/p}\left([-r,0]\right) \text{ with } p\in(1,+\infty]$$

and thus, we obtain in any case

$$\sup\left\{ \|x_t\|_X : t\in[0,r], x_0\in X, \|x_0\|_X \le R \right\} < +\infty \tag{23}$$

In order to establish (22), let arbitrary $R>0$ be given and define $B:=\left\{ x\in X : \|x\|_X \le R \right\}$ (the ball of radius $R>0$ in $X$). Consider next the closure $\bar{B}$ of $B$ in $C^0\left([-r,0];\Re^n\right)$. Theorem 8.8 on pages 212-213 in [2] implies that $\bar{B}$ is compact in the topology of $C^0\left([-r,0];\Re^n\right)$ and $\bar{B}\subseteq C^{0,1-1/p}\left([-r,0]\right)$.

By virtue of continuity of solutions with respect to the initial conditions for each $x\in C^{0,1-1/p}\left([-r,0]\right)$ there exists $\eta(x)>0$ such that $\|\phi(t,y)-\phi(t,x)\|_\infty \le 1$ for all $t\in[0,r]$ and for all $y\in C^0\left([-r,0];\Re^n\right)$ with $\|y-x\|_\infty < \eta(x)$. For each $x\in C^{0,1-1/p}\left([-r,0]\right)$ the set $N(x):=\left\{ y\in C^0\left([-r,0];\Re^n\right) : \|y-x\|_\infty < \eta(x) \right\}$ is an open neighborhood of $x$ in the topology of $C^0\left([-r,0];\Re^n\right)$. Consequently, the sets $N(x)$ for $x\in\bar{B}$ constitute an open cover of $\bar{B}$. Therefore, by compactness of the set $\bar{B}$ there exists a (finite) positive integer $m$ and points $x_i\in\bar{B}$, $i=1,...,m$ such that $\bar{B}\subseteq\bigcup_{i=1}^{m} N(x_i)$.

Define

$$\bar{\rho}:=\max_{i=1,...,m}\left(\sup\left\{ \|\phi(t,x_i)\|_\infty : t\in[0,r] \right\}\right)+1 \tag{24}$$

Let arbitrary $y\in B$ be given. Since $B\subseteq\bar{B}\subseteq\bigcup_{i=1}^{m} N(x_i)$, there exists $i\in\{1,...,m\}$ such that $y\in N(x_i)=\left\{ z\in C^0\left([-r,0];\Re^n\right) : \|z-x_i\|_\infty < \eta(x_i) \right\}$. Consequently, by virtue of definition (24), the fact that $\|\phi(t,y)-\phi(t,x_i)\|_\infty \le 1$ for all $t\in[0,r]$ and the triangle inequality, we get:

$$\begin{aligned}
&\sup\left\{ \|\phi(t,y)\|_\infty : t\in[0,r] \right\} \\
&\le \sup\left\{ \|\phi(t,y)-\phi(t,x_i)\|_\infty + \|\phi(t,x_i)\|_\infty : t\in[0,r] \right\} \\
&\le \sup\left\{ \|\phi(t,y)-\phi(t,x_i)\|_\infty : t\in[0,r] \right\} + \sup\left\{ \|\phi(t,x_i)\|_\infty : t\in[0,r] \right\} \\
&\le 1+\sup\left\{ \|\phi(t,x_i)\|_\infty : t\in[0,r] \right\} \le \bar{\rho}
\end{aligned} \tag{25}$$

It follows from (25) and the fact that $y\in B$ is arbitrary that the following inequality holds:

$$\sup\left\{ \|\phi(t,y)\|_\infty : t\in[0,r], y\in B \right\} \le \bar{\rho} < +\infty \tag{26}$$

Since $B:=\left\{ x\in X : \|x\|_X \le R \right\}$, it follows that (26) implies (22). The proof is complete. $\lhd$

**Proof of Theorem 3:** Lemma 2.1 on page 58 in [13] shows the equivalence

(UGA-X) $\wedge$ (RFC-X) $\Leftrightarrow$ (UGAS-X)

Since Property (RFC-X) holds (Theorem 2), the implication (LS-X) $\wedge$ (GA-X) $\Rightarrow$ (UGAS-X) follows directly from the implication (LS-X) $\wedge$ (GA-X) $\Rightarrow$ (UGA-X).

Thus, it suffices to show the implication (LS-X) $\wedge$ (GA-X) $\Rightarrow$ (UGA-X).

Let arbitrary $\varepsilon, \rho > 0$ be given.

By virtue of Lyapunov Stability (Property (LS-X)), there exists $\delta > 0$ so that

$$\sup\left\{ \|\phi(t,x_0)\|_X : t \geq 0, x_0 \in X, \|x_0\|_X \leq \delta \right\} \leq \varepsilon \tag{27}$$

By virtue of property (RFC-X) (Theorem 2) there exists $R > 0$ so that

$$\sup\left\{ \|\phi(t,x_0)\|_X : t \in [0,r], x_0 \in X, \|x_0\|_X \leq \rho \right\} \leq R \tag{28}$$

Let $M := \max\left\{ |f(x)| : \|x\|_\infty \leq R \right\}$ and define the set

$$B := \left\{ x \in W^{1,\infty}\left([-r,0];\Re^n\right) : \|\dot{x}\|_\infty \leq M, \|x\|_\infty \leq R \right\} \tag{29}$$

Then by virtue of the Arzela-Ascoli theorem and the fact that a function $x:[-r,0] \to \Re^n$ is in $B$ if and only if $x$ is Lipschitz with Lipschitz constant $M$ and $\|x\|_\infty \leq R$, the set $B$ is compact in the topology of $C^0\left([-r,0];\Re^n\right)$.

Definition (29), estimate (28) and (1) imply that for every $x_0 \in X$ with $\|x_0\|_X \leq \rho$ it holds that $\phi(r,x_0) \in B$. Therefore, by virtue of the semigroup property, in order to prove that there exists $\bar{T} > 0$ such that

$$\sup\left\{ \|\phi(t,x_0)\|_X : t \geq \bar{T}, x_0 \in X, \|x_0\|_X \leq \rho \right\} \leq \varepsilon \tag{30}$$

it suffices to show that there exists $T > 0$ such that

$$\sup\left\{ \|\phi(t,x_0)\|_X : t \geq T, x_0 \in B \right\} \leq \varepsilon . \tag{31}$$

Inequality (30) follows from inequality (31) by setting $\bar{T} = T + r$.

Let arbitrary $x_0 \in X$ be given. By virtue of Global Attractivity (Property (GA-X)), there exists $\tilde{T}(x_0) > 0$ such that $\|\phi(t,x_0)\|_X \leq \dfrac{\delta}{2\left(1 + r^{1/p} L(\varepsilon)\right)}$ for all $t \geq \tilde{T}(x_0)$, where $L$ is the function involved in (2).

By virtue of continuity of solutions with respect to the initial conditions for each $x_0 \in X$ there exists $\eta(x_0) > 0$ such that $\|\phi(t,y) - \phi(t,x_0)\|_\infty \leq \dfrac{\delta}{2\left(1 + r^{1/p} L(\varepsilon)\right)}$ for all $t \in \left[0, \tilde{T}(x_0) + r\right]$ and for all

$y \in C^0\left([-r,0];\Re^n\right)$ with $\|y-x_0\|_\infty < \eta(x_0)$. Thus, by virtue of the triangle inequality and the fact that $\|\phi(t,x_0)\|_\infty \le \|\phi(t,x_0)\|_W$ we get for all $y \in C^0\left([-r,0];\Re^n\right)$ with $\|y-x_0\|_\infty < \eta(x_0)$ and for all $t \in \left[\tilde{T}(x_0), \tilde{T}(x_0)+r\right]$:

$$\begin{aligned} \left\|\phi(t,y)\right\|_\infty &\le \left\|\phi(t,y)-\phi(t,x_0)\right\|_\infty + \left\|\phi(t,x_0)\right\|_\infty \\ &\le \frac{\delta}{2\left(1+r^{1/p}L(\varepsilon)\right)} + \frac{\delta}{2\left(1+r^{1/p}L(\varepsilon)\right)} \le \frac{\delta}{1+r^{1/p}L(\varepsilon)} \end{aligned} \tag{32}$$

It follows from (1), (2) and the fact that $L:\Re_+ \to \Re_+$ is non-decreasing that the following inequality holds for all $t \in \left[r, \tilde{T}(x_0)+r\right]$:

$$\|\phi(t,y)\|_X \le \left(1+r^{1/p}L\left(\max_{t-r\le s\le t}\left(\|\phi(s,y)\|_\infty\right)\right)\right)\max_{t-r\le s\le t}\left(\|\phi(s,y)\|_\infty\right) \tag{33}$$

Therefore, for each $x_0 \in X$ and for all $y \in C^0\left([-r,0];\Re^n\right)$ with $\|y-x_0\|_\infty < \eta(x_0)$ we obtain from (32) and (33) (using again the fact that $L:\Re_+ \to \Re_+$ is non-decreasing and the fact that $\delta \le \varepsilon$ which is a consequence of (27)):

$$\left\|\phi(\tilde{T}(x_0)+r,y)\right\|_X \le \frac{\left(1+r^{1/p}L(\delta)\right)\delta}{1+r^{1/p}L(\varepsilon)} \le \delta \tag{34}$$

Consequently, by virtue of (27), (34) and the semigroup property, we get for each $x_0 \in X$:

$$\begin{aligned} &\sup\left\{ \|\phi(t,y)\|_X : t \ge \tilde{T}(x_0)+r,\, y \in C^0\left([-r,0];\Re^n\right), \|y-x_0\|_\infty < \eta(x_0) \right\} \\ &= \sup\left\{ \left\|\phi\left(t-\tilde{T}(x_0)-r, \phi\left(\tilde{T}(x_0)+r,y\right)\right)\right\|_X : t \ge \tilde{T}(x_0)+r,\, y \in C^0\left([-r,0];\Re^n\right), \|y-x_0\|_\infty < \eta(x_0) \right\} \\ &\le \sup\left\{ \left\|\phi(s,z)\right\|_X : s \ge 0,\, z \in X,\, \|z\|_X \le \delta \right\} \le \varepsilon \end{aligned} \tag{35}$$

For each $x \in X$ the set $N(x) := \left\{ y \in C^0\left([-r,0];\Re^n\right) : \|y-x\|_\infty < \eta(x) \right\}$ is an open neighborhood of $x$ in the topology of $C^0\left([-r,0];\Re^n\right)$. Consequently, the sets $N(x)$ for $x \in B$ constitute an open cover of $B$. Therefore, by compactness of the set $B$ there exists a (finite) positive integer $m$ and points $x_i \in B$, $i = 1,...,m$ such that $B \subseteq \bigcup_{i=1}^{m} N(x_i)$.

Define

$$T := r + \max_{i=1,...,m}\left(\tilde{T}(x_i)\right) \tag{36}$$

Let arbitrary $y \in B$ be given. Since $B \subseteq \bigcup_{i=1}^{m} N(x_i)$, there exists $i \in \{1,...,m\}$ such that $y \in N(x_i) \cap B = \left\{ z \in B : \|z - x_i\|_\infty < \eta(x_i) \right\}$. Consequently, by virtue of definition (36) and inequality (35), we get:

$$\sup\left\{ \|\phi(t,y)\|_X : t \geq T \right\} \leq \sup\left\{ \|\phi(t,z)\|_X : t \geq r + \tilde{T}(x_i),\, z \in B,\, \|z - x_i\|_\infty < \eta(x_i) \right\} \leq \varepsilon \tag{37}$$

It follows from (37) and the fact that $y \in B$ is arbitrary, that inequality (31) holds. The proof is complete. $\lhd$

**Proof of Theorem 4:** Since the implications (UGAS) $\Rightarrow$ (Q-X) and (UGAS-X) $\Rightarrow$ (Q-X) are trivial, we focus on the proof of the implication (Q-X) $\Rightarrow$(UGAS-X).

Suppose that there exists $\sigma \in KL$ such that the following estimate

$$\|x_t\|_\infty \leq \sigma\left( \|x_0\|_X , t \right), \text{ for all } t \geq 0,\ x_0 \in X . \tag{38}$$

Combining (1), (2) and (38) we obtain that the estimate

$$|\dot{x}(t)| = |f(x_t)| \leq L\left(\sigma\left( \|x_0\|_X , 0 \right)\right)\|x_t\|_\infty \tag{39}$$

holds for all $t \geq 0$, $x_0 \in X$ .

We next distinguish the following cases.

<u>Case 1:</u> $t \in [0,r]$. For this case, it follows by virtue of (38), (39), the facts that $\|\dot{x}_0\|_p \leq \|x_0\|_X$ for all $x_0 \in X$ when $X = W^{1,p}\left([-r,0]\right)$ and $\sup_{l,s\in[-r,0],l\neq s}\left( \frac{|x_0(l) - x_0(s)|}{|l-s|^{1-1/p}} \right) \leq \|x_0\|_X$ for all $x_0 \in X$ when $X = C^{0,1-1/p}\left([-r,0]\right)$ and the fact that $\sigma(\ s, 0\ ) \geq s$ for all $s \geq 0$ that the estimates

$$\begin{aligned}
\|\dot{x}_t\|_p &\leq \left(1 + r^{1/p}\right)\max\left( \|\dot{x}_0\|_p , \sup_{0<s<t}\left(|\dot{x}(s)|\right) \right) \\
&\leq \left(1 + r^{1/p}\right)\max\left( \|x_0\|_X , L\left(\sigma\left( \|x_0\|_X , 0 \right)\right)\sigma\left( \|x_0\|_X , 0 \right) \right) \\
&\leq \left(1 + r^{1/p}\right)\exp\left(r - t\right)\max\left(1, L\left(\sigma\left( \|x_0\|_X , 0 \right)\right)\right)\sigma\left( \|x_0\|_X , 0 \right)
\end{aligned}$$

when $X = W^{1,p}\left([-r,0]\right)$

$$\sup_{l,s\in[-r,0],l\neq s}\left(\frac{|x_t(l)-x_t(s)|}{|l-s|^{1-1/p}}\right)=\sup_{l,s\in[-r,0],l\neq s}\left(\frac{|x(t+l)-x(t+s)|}{|l-s|^{1-1/p}}\right)$$
$$\leq \sup_{l,s\in[-r,0],l\neq s}\left(\frac{|x_0(l)-x_0(s)|}{|l-s|^{1-1/p}}\right)+r^{1/p}\sup_{0<s<t}\left(|\dot{x}(s)|\right)$$
$$\leq \|x_0\|_X + r^{1/p}L\left(\sigma\left(\|x_0\|_X,0\right)\right)\sigma\left(\|x_0\|_X,0\right)$$
$$\leq \left(1+r^{1/p}\right)\exp\left(r-t\right)\max\left(1,L\left(\sigma\left(\|x_0\|_X,0\right)\right)\right)\sigma\left(\|x_0\|_X,0\right)$$
$$\text{when } X=C^{0,1-1/p}\left([-r,0]\right)$$

hold for all $t\in(0,r)$ and all $x_0\in X$. We notice next that inequalities

$$\|\dot{x}_t\|_p \leq \left(1+r^{1/p}\right)\exp\left(r-t\right)\max\left(1,L\left(\sigma\left(\|x_0\|_X,0\right)\right)\right)\sigma\left(\|x_0\|_X,0\right)$$
$$\text{when } X=W^{1,p}\left([-r,0]\right) \tag{40}$$

$$\sup_{l,s\in[-r,0],l\neq s}\left(\frac{|x_t(l)-x_t(s)|}{|l-s|^{1-1/p}}\right)\leq \left(1+r^{1/p}\right)\exp\left(r-t\right)\max\left(1,L\left(\sigma\left(\|x_0\|_X,0\right)\right)\right)\sigma\left(\|x_0\|_X,0\right)$$
$$\text{when } X=C^{0,1-1/p}\left([-r,0]\right) \tag{41}$$

hold for all $t\in[0,r]$ and all $x_0\in X$.

<u>Case 2:</u> $t>r$. For this case, it follows by virtue of (38), (39) that the estimates

$$\|\dot{x}_t\|_p \leq r^{1/p}\max_{t-r\leq s\leq t}\left(|\dot{x}(s)|\right)\leq r^{1/p}L\left(\sigma\left(\|x_0\|_X,0\right)\right)\max_{t-r\leq s\leq t}\left(\|x_s\|\right)$$
$$\leq r^{1/p}L\left(\sigma\left(\|x_0\|_X,0\right)\right)\sigma\left(\|x_0\|_X,t-r\right)$$
$$\leq \left(1+r^{1/p}\right)\max\left(1,L\left(\sigma\left(\|x_0\|_X,0\right)\right)\right)\sigma\left(\|x_0\|_X,t-r\right)$$
$$\text{when } X=W^{1,p}\left([-r,0]\right) \tag{42}$$

$$\sup_{l,s\in[-r,0],l\neq s}\left(\frac{|x_t(l)-x_t(s)|}{|l-s|^{1-1/p}}\right)=\sup_{l,s\in[-r,0],l\neq s}\left(\frac{|x(t+l)-x(t+s)|}{|l-s|^{1-1/p}}\right)$$
$$\leq r^{1/p}\max_{t-r\leq s\leq t}\left(|\dot{x}(s)|\right)\leq r^{1/p}L\left(\sigma\left(\|x_0\|_X,0\right)\right)\max_{t-r\leq s\leq t}\left(\|x_s\|\right)$$
$$\leq r^{1/p}L\left(\sigma\left(\|x_0\|_X,0\right)\right)\sigma\left(\|x_0\|_X,t-r\right)$$
$$\leq \left(1+r^{1/p}\right)\max\left(1,L\left(\sigma\left(\|x_0\|_X,0\right)\right)\right)\sigma\left(\|x_0\|_X,t-r\right)$$
$$\text{when } X=C^{0,1-1/p}\left([-r,0]\right) \tag{43}$$

hold for all $t>r$ and all $x_0\in X$.

Define the $KL$ function

$$\omega(s,t) := \sigma(s,t) + \left(1+r^{1/p}\right)\max\left(1, L\left(\sigma(s,0)\right)\right)\begin{cases}\exp(r-t)\sigma(s,0) & \text{for } t\in[0,r] \\ \sigma(s,t-r) & \text{for } t>r\end{cases}$$

$$\text{for all } t,s\geq 0 . \tag{44}$$

Combining estimates (38), (40), (41), (42) and (43), we conclude that the estimate $\|x_t\|_X \leq \omega\left(\|x_0\|_X, t\right)$ holds for all $t\geq 0$ and all $x_0\in X$. The proof is complete. $\lhd$

**Proof of Theorem 5:** Suppose that there exist a functional $V: X\to\Re_+$ which is Lipschitz on bounded sets of $X$ and functions $a_1,a_2\in K_\infty$ such that inequalities (3), (4) hold for all $x\in X$. Combining (3), (4) we get the estimate $\|\phi(t,x)\|_X \leq a_1^{-1}\left(\exp(-t)a_2\left(\|x\|_X\right)\right)$ for all $t\geq 0$, $x\in X$. It is therefore clear that Property (UGAS-X) holds with $\sigma(s,t) := a_1^{-1}\left(\exp(-t)a_2(s)\right)$ for all $t,s\geq 0$.

Next suppose that Property (UGAS-X) holds. In order to show that there exist a functional $V: X\to\Re_+$ which is Lipschitz on bounded sets of $X$ and functions $a_1,a_2\in K_\infty$ such that inequalities (3), (4) hold for all $x\in X$ it suffices to show that the following property holds.

**(H)** *For every $T,R>0$ there exists $M>0$ such that the estimate $\|x_t-y_t\|_X \leq M\|x_0-y_0\|_X$ holds for all $t\in[0,T]$, $x_0,y_0\in X$ with $\|x_0\|_X\leq R$, $\|y_0\|_X\leq R$, where $x(t),y(t)$ are the solutions of (1) with initial conditions $x_0,y_0\in X$, respectively.*

Indeed, Property (H) is equivalent to Property (REG2) on page 130 in [13]. The existence of a functional $V: X\to\Re_+$ which is Lipschitz on bounded sets of $X$ and functions $a_1,a_2\in K_\infty$ such that inequalities (3), (4) hold for all $x\in X$ is a direct consequence of Theorem 3.4 on pages 135-136 in [13].

Let $T,R>0$ be given (arbitrary). Let $x_0,y_0\in X$ with $\|x_0\|_X\leq R$, $\|y_0\|_X\leq R$ be given (arbitrary). Let $x(t),y(t)$ be the solutions of (1) with initial conditions $x_0,y_0\in X$, respectively. Using (1) we get for all $t\in[0,T]$:

$$|x(t)-y(t)| = \left|x(0)-y(0)+\int_0^t \left(f(x_s)-f(y_s)\right)ds\right| \tag{45}$$

Using the triangle inequality, (1), (2) and the facts that $\|x_0\|_X\leq R$, $\|y_0\|_X\leq R$, $\|x_t\|_\infty \leq \|x_t\|_X \leq \sigma\left(\|x_0\|_X, t\right) \leq \sigma\left(\|x_0\|_X, 0\right)$, $\|y_t\|_\infty \leq \|y_t\|_X \leq \sigma\left(\|y_0\|_X, t\right) \leq \sigma\left(\|y_0\|_X, 0\right)$ for all $t\geq 0$ (consequences of Property (UGAS-X)), we obtain from (45) for all $t\in[0,T]$:

$$|x(t)-y(t)| \leq |x(0)-y(0)| + L\left(\sigma(R,0)\right)\int_0^t \|x_s-y_s\|_\infty\, ds \tag{46}$$

$$|\dot{x}(t)-\dot{y}(t)| \leq L\left(\sigma(R,0)\right)\|x_t-y_t\|_\infty \tag{47}$$

Estimate (46) implies the following estimate for all $t\in[0,T]$:

$$\|x_t-y_t\|_\infty \leq \|x_0-y_0\|_\infty + L\left(\sigma(R,0)\right)\int_0^t \|x_s-y_s\|_\infty\, ds \tag{48}$$

Using the Gronwall-Bellman Lemma and (48) we get for all $t \in [0,T]$:

$$\|x_t - y_t\|_\infty \le \exp\left(L\left(\sigma(R,0)\right)T\right)\|x_0 - y_0\|_\infty \tag{49}$$

We also get from (47) and (49) for all $t \in [0,T]$:

$$|\dot{x}(t) - \dot{y}(t)| \le L\left(\sigma(R,0)\right)\exp\left(L\left(\sigma(R,0)\right)T\right)\|x_0 - y_0\|_\infty \tag{50}$$

Consequently, we get from (50) for all $t \in [0,T]$:

$$\|\dot{x}_t - \dot{y}_t\|_p \le \left(1 + r^{1/p}\right)\max\left(\|\dot{x}_0 - \dot{y}_0\|_p \,, L\left(\sigma(R,0)\right)\exp\left(L\left(\sigma(R,0)\right)T\right)\|x_0 - y_0\|_\infty\right)$$
$$\text{when } X = W^{1,p}\left([-r,0]\right) \tag{51}$$

$$\sup_{l,s\in[-r,0],l\neq s}\left(\frac{|x_t(l) - y_t(l) - x_t(s) + y_t(s)|}{|l-s|^{1-1/p}}\right)$$
$$\le \sup_{l,s\in[-r,0],l\neq s}\left(\frac{|x_0(l) - y_0(l) - x_0(s) + y_0(s)|}{|l-s|^{1-1/p}}\right) + r^{1/p}L\left(\sigma(R,0)\right)\exp\left(L\left(\sigma(R,0)\right)T\right)\|x_0 - y_0\|_\infty$$
$$\text{when } X = C^{0,1-1/p}\left([-r,0]\right) \tag{52}$$

Using (49), (51), (52) and the facts that

$\|x_0 - y_0\|_\infty \le \|x_0 - y_0\|_X$ for every considered state space,

$\|\dot{x}_0 - \dot{y}_0\|_p \le \|x_0 - y_0\|_X$ and $\|x_t - y_t\|_X = \|x_t - y_t\|_\infty + \|\dot{x}_t - \dot{y}_t\|_p$ when $X = W^{1,p}\left([-r,0]\right)$,

$\sup_{l,s\in[-r,0],l\neq s}\left(\frac{|x_0(l) - y_0(l) - x_0(s) + y_0(s)|}{|l-s|^{1-1/p}}\right) \le \|x_0 - y_0\|_X$ and

$\|x_t - y_t\|_X = \max\left(\|x_t - y_t\|_\infty \,, \sup_{t,s\in[-r,0],t\neq s}\left(\frac{|x_t(l) - y_t(l) - x_t(s) + y_t(s)|}{|t-s|^a}\right)\right)$ when $X = C^{0,1-1/p}\left([-r,0]\right)$,

we obtain for all $t \in [0,T]$:

$$\|x_t - y_t\|_X \le \left(1 + \left(1 + r^{1/p}\right)\max\left(1, L\left(\sigma(R,0)\right)\exp\left(L\left(\sigma(R,0)\right)T\right)\right)\right)\|x_0 - y_0\|_X \tag{53}$$

Therefore, Property (H) holds with $M = 1 + \left(1 + r^{1/p}\right)\max\left(1, L\left(\sigma(R,0)\right)\exp\left(L\left(\sigma(R,0)\right)T\right)\right)$. The proof is complete. $\triangleleft$

**Proof of Lemma 1:** Consider the solution $x(t)$ of (1) with initial condition $x_0 \in X$ and let $t_{\max} \in (0,+\infty]$ denote its maximal time of existence. Given $T \in (0, t_{\max})$, the solution $x(t)$ exists for $t \in [0,T]$. Then we may define $R := \max_{0\le t\le T}\left(|\dot{x}(t)|\right)$ (the fact that $\max_{0\le t\le T}\left(|\dot{x}(t)|\right) < +\infty$ follows from (2) and continuity of the mapping $t \mapsto x_t$ in the topology of $C^0\left([-r,0];\Re^n\right)$; see Lemma 2.1 on page 40 in [9]). Then it holds that

$$\left( \int_{-r}^{T} |\dot{x}(s)|^p \, ds \right)^{1/p} = \left( \int_{-r}^{0} |\dot{x}(s)|^p \, ds + \int_{0}^{T} |\dot{x}(s)|^p \, ds \right)^{1/p} \leq \left( \|\dot{x}_0\|_p^p + TR^p \right)^{1/p} \tag{54}$$

Thus $\dot{x} \in L^p\left((-r,T);\Re^n\right)$. Let arbitrary $\varepsilon > 0$ be given. Using density of $C^0\left([-r,T];\Re^n\right)$ in $L^p\left((-r,T);\Re^n\right)$ there exists $h \in C^0\left([-r,T];\Re^n\right)$ such that $\left( \int_{-r}^{T} |h(s) - \dot{x}(s)|^p \, ds \right)^{1/p} \leq \frac{\varepsilon}{4}$. Let arbitrary $t_1, t_2 \in [0,T]$ with $t_2 \geq t_1$ be given. We have

$$\begin{aligned}
\left\| \dot{x}_{t_2} - \dot{x}_{t_1} \right\|_p &= \left( \int_{-r}^{0} |\dot{x}(t_2 + s) - \dot{x}(t_1 + s)|^p \, ds \right) \\
&\leq \left( \int_{t_1 - r}^{t_1} |\dot{x}(t_2 - t_1 + s) - h(t_2 - t_1 + s)|^p \, ds \right)^{1/p} \\
&+ \left( \int_{t_1 - r}^{t_1} |h(t_2 - t_1 + s) - h(s)|^p \, ds \right)^{1/p} + \left( \int_{t_1 - r}^{t_1} |h(s) - \dot{x}(s)|^p \, ds \right)^{1/p} \\
&\leq 2 \left( \int_{-r}^{T} |h(s) - \dot{x}(s)|^p \, ds \right)^{1/p} + \left( \int_{t_1 - r}^{t_1} |h(t_2 - t_1 + s) - h(s)|^p \, ds \right)^{1/p} \\
&\leq \frac{\varepsilon}{2} + \left( \int_{t_1 - r}^{t_1} |h(t_2 - t_1 + s) - h(s)|^p \, ds \right)^{1/p}
\end{aligned}$$

Moreover, by virtue of continuity of $h \in C^0\left([-r,T];\Re^n\right)$, there exists $\delta > 0$ such that $|h(s+\tau) - h(s)| \leq \frac{\varepsilon}{4r^{1/p}}$ for all $s \in [-r,T]$, $\tau \in [0,\delta]$ with $s + \tau \leq T$ (uniform continuity on compact sets). It follows from the previous estimates that $\left\| \dot{x}_{t_2} - \dot{x}_{t_1} \right\|_p \leq 3\varepsilon / 4$ for all $t_1, t_2 \in [0,T]$ with $t_2 \in [t_1, t_1 + \delta]$. By continuity of the mapping $t \mapsto x_t$ in the topology of $C^0\left([-r,0];\Re^n\right)$ (see Lemma 2.1 on page 40 in [9]) there exists $\bar{\delta} > 0$ such that $\left\| x_{t_2} - x_{t_1} \right\|_\infty \leq \frac{\varepsilon}{4}$ for all $t_1, t_2 \in [0,T]$ with $t_2 \in [t_1, t_1 + \bar{\delta}]$. Therefore, we get $\left\| x_{t_2} - x_{t_1} \right\|_X \leq \varepsilon$ for all $t_1, t_2 \in [0,T]$ with $0 \leq t_2 - t_1 \leq \min\left(\delta, \bar{\delta}\right)$. Since $T \in (0, t_{\max})$ and $\varepsilon > 0$ are arbitrary we conclude that the mapping $t \mapsto x_t = \phi(t, x_0)$ is continuous in the topology of $X = W^{1,p}\left([-r,0]\right)$. The proof is complete. $\lhd$

**Proof of Theorem 6:** Implication (i) $\Rightarrow$ (ii) is a consequence of Theorem 3, implication (ii) $\Rightarrow$ (iii) is a consequence of Theorem 5 and implication (iii) $\Rightarrow$ (iv) is trivial. Therefore, it suffices to prove implication (iv) $\Rightarrow$ (i). Let arbitrary $x_0 \in X$ and consider the solution $x(t)$ of (1) with initial condition $x_0 \in X$ defined for $t \in [0, t_{\max})$ with $t_{\max} \in (0, +\infty]$ being the maximal existence time of the solution. It is known that if $t_{\max} < +\infty$ then $\limsup_{t \to t_{\max}^-} \left( |x(t)| \right) = +\infty$.

Using Lemma 1, and the fact that $V$ is Lipschitz on bounded sets of $X$, we get that the mapping $t \mapsto V(x_t)$ is continuous on $[0, t_{\max})$. Using (10) and the semigroup property we get that $\limsup_{h \to 0^+} \left( h^{-1} \left( V(x_{t+h}) - V(x_t) \right) \right) \leq 0$ for all $t \in [0, t_{\max})$. Lemma 2.12 on pages 77-78 in [13] implies that

$$V(x_t) \leq V(x_0), \text{ for all } t \in [0, t_{\max}). \tag{55}$$

Combining (55) with (9) gives the estimate

$$|x(t)| \leq a_1^{-1} \left( a_2 \left( \|x_0\|_X \right) \right), \text{ for all } t \in [0, t_{\max}). \tag{56}$$

It follows from (56) that $t_{\max} = +\infty$, meaning that the system is forward complete. Moreover, estimate (56) in conjunction with (2) shows that the property (LS-X) holds.

Using (10) and the semigroup property we get that $\limsup_{h \to 0^+} \left( h^{-1} \left( V(x_{t+h}) - V(x_t) \right) \right) \leq -Q(x(t))$ for all $t \geq 0$. By virtue of a slight variation of Theorem 3 in [8] (considering continuous functions with right upper Dini derivative bounded from above on their domains), we get for all $t \geq 0$:

$$V(x_t) + \int_0^t Q(x(s)) ds \leq V(x_0) \tag{57}$$

It follows from continuity of $Q$, estimates (56), (57), the fact that $Q$ is positive definite, the fact that $|\dot{x}(t)|$ is bounded for $t \geq 0$ (a consequence of (56)) and Barbălat's Lemma (see [14]) that $\lim_{t \to +\infty} \left( x(t) \right) = 0$. Hence, it follows from (2) that $\lim_{t \to +\infty} \left( \|x_t\|_X \right) = 0$. Thus, the property (GA-X) holds. The proof is complete. $\triangleleft$

## 4. Conclusions

In this paper, we have shown that GAS and UGAS are equivalent properties for time-invariant delay systems provided that the RFC property holds. To our knowledge, it has not yet been proved or disproved that RFC automatically holds under forward completeness for this class of systems (as it does in finite dimension). A positive answer to this question would show that GAS and UGAS are equivalent properties also for time-delay systems. We therefore believe this is a central question that deserves deeper investigations and refer the reader to [20] for further discussions on this matter.

When the state space is the Sobolev space $W^{1,p}\left([-r,0]\right)$ with $p \in (1, +\infty]$, we have also shown that under a forward completeness property, RFC does hold and GAS and UGAS are equivalent. When the state space is the Hölder space $C^{0,1-1/p}\left([-r,0]\right)$ with $p \in (1, +\infty]$, we have also shown that under the usual forward completeness property, RFC does hold and GAS and UGAS are equivalent, just like in the finite-dimensional case. Finally, in the Sobolev space $W^{1,p}\left([-r,0]; \Re^n\right)$ with $p \in (1, +\infty)$, we have provided a Lyapunov-Krasovskii characterization of GAS.